\documentclass{amsart}
\usepackage{amssymb}
\usepackage{amsmath}
\usepackage{amscd}
\usepackage{bbm}

\newcommand\Rs{{\mathbb R}}

\newcommand\supp{\mathop {\fam 0 supp}\nolimits}
\newcommand\dd{\partial}

\newcommand\ddiv{\mathop {\fam 0 div}\nolimits}

\newcommand\CC{{\mathcal C}}

\begin{document}

\title{A  steady Euler flow with compact support}

\author{A.~V.~Gavrilov}

             \email{gavrilov19@gmail.com}           


\maketitle

\begin{abstract}
A nontrivial  smooth steady incompressible Euler flow in three dimensions  with compact support is constructed.
Another   uncommon property of this solution  is  the dependence between the Bernoulli function and  the  pressure.
\end{abstract}

\section{Introduction}

A steady flow of an ideal  fluid in $\Rs^3$  is a solution of  the Euler equation
$$(u\cdot\nabla) u=-\nabla p,\,\ddiv u=0.$$
At present, it is not known if  smooth nonzero   solutions of this equation with compact support  $0\neq u\in C_0^\infty(\Rs^3)$ exist  \cite{N}\cite{NV}. 
The  problem is trivial in two dimensions where there are obvious vortex-like solutions. In three dimensions, 
only a few results are known, all on the negative side.   It is known that such a flow cannot be Beltrami 
\cite{N}\cite{CC}, and it  cannot be  axisymmetric  without swirl  \cite{JX}. Recently, 
Nadirashvili and Vladut found some other  restrictions  \cite{NV}.

Apparently\footnote{The author is no expert in this area.},  weak  solutions with compact support may be constructed using methods of  \cite{CS}. 
Also, there is a considerable literature about  vortex rings  which are solutions with compactly supported vorticity (e.g. \cite{AS}).
Opinions have been expressed that in three dimensions  there are no smooth solutions with compact support besides $u=0$.
The main goal of this paper is to show that it is not true.

{\bf Theorem} \emph{There exists a nontrivial smooth steady Euler flow  in $\Rs^3$
with support in an arbitrarily  small neighbourhood of a circle.}

We give below  an explicit description  of   an axisymmetric    flow with  compact support. 
This solution  has also  other unusual  properties discussed in the last section.

\emph{Acknowledgement.} The author would like to thank the  anonymous referee for  
pointing  out  the  interesting   work of  Khesin, Kuksin, and  Peralta-Salas  \cite{KKP} as well as some properties of the given solution.

\section{Some differential equations}

In this section we find solutions of some  differential equations  which will be  used later.

\subsection{Ordinary differential equation}

$$$$
{\bf Lemma 1}

The  singular  Cauchy problem
\begin{equation}
3x\psi''+6x(\psi')^3-4\psi(\psi')^2-3\psi'=0;\,\psi(0)=1, \psi'(0)=-\frac{3}{4} \label{e1}
\end{equation}
has  a unique analytic  solution $\psi(x)$ in a neighbourhood of $x=0$.

\begin{proof}
The equation  \eqref{e1} becomes  first order   in variables 
$t=x\psi^{-2}$ and $v=\psi\psi',$
$$t\frac{dv}{dt}=v\left(\frac{4}{3}v+1\right)+\frac{tv^2(2v+9)}{3(1-2tv)}.$$
Denoting $w=v+\frac{3}{4}$  we may write this equation as 
$t\frac{dw}{dt}=-w+f(t,w)$
where the function $f$ is analytic and $f(0,0)=\frac{\dd}{\dd w}f(0,0)=0$.  By  \cite[Theorem 11.1]{H}  there is  an unique analytic   solution $v(t)$ such that $v(0)=-\frac{3}{4}$. The   Cauchy problem $\frac{d\psi}{d x}=\frac{1}{\psi}v\left(\frac{x}{\psi^2}\right),\,\psi(0)=1$  clearly  has  an unique solution.
\end{proof}

From now   $\psi$  always means  the function  defined by \eqref{e1}.   The Taylor series of this function is
$$\psi(x)=1-\frac{3}{4}x+\frac{9}{128}x^2-\frac{21}{1024}x^3+\frac{1035}{131072}x^4-\frac{1809}{524288}x^5+O(x^6).$$

 Denote  
$$F(x, \alpha)=-2x\psi(\alpha)+2x^3,\,\,H(\alpha)=6\alpha\left(\frac{1}{\psi'(\alpha)}+2\psi(\alpha)\right),$$  
$$G(x, \alpha)=12x^2 \alpha - F^2(x,\alpha)-H(\alpha).$$
Note that at $(x, \alpha)=(1,0)$ we have $F=G=0$ and
$$\frac{\dd F}{\dd x}=4,\,\frac{\dd F}{\dd\alpha}=\frac{3}{2},\,\frac{\dd G}{\dd x}=0,\,\frac{\dd G}{\dd\alpha}=8.$$

We will also need  the following fact  
$$$$
{\bf Lemma 2}
The  functions $F,G$ satisfy
\begin{equation}
\frac{\dd G}{\dd x}+F\frac{\dd G}{\dd\alpha}=2G\frac{\dd F}{\dd\alpha},  \label{e2a}
\end{equation}

\begin{equation}
x\frac{\dd F}{\dd x}-F=4x^3.  \label{e2b}
\end{equation}

\begin{proof}

The  part \eqref{e2b} is trivial; \eqref{e2a}  boils down to the formula 
$$H'(\alpha)=24\alpha\psi'(\alpha)+4\psi(\alpha),$$
equivalent to \eqref{e1}. 
\end{proof}

\subsection{Partial differential equations}

$$$$
{\bf Lemma 3}
The system\footnote{Note that  $\frac{\dd }{\dd x}$ means $\frac{\dd }{\dd x}|_{\alpha}$ when applied to $F$ or $G$ but  $\frac{\dd }{\dd x}|_{y}$  when applied to $\alpha$. To avoid a (very common)  inconsistency  in notation we write this partial derivative as $\frac{\dd f}{\dd x}$ for the former and $\frac{\dd }{\dd x}f$ for the latter.}  
\begin{equation}
\frac{\dd }{\dd x}\alpha =F(x, \alpha),\,
\left(\tfrac{\dd }{\dd y} \alpha\right)^2=G(x, \alpha),
\label{e2}
\end{equation}

has a unique analytic solution $\alpha(x,y)$  in a neighbourhood of the point $(x,y)=(1,0)$ 
such that   $\alpha(1,0)=0$ and  $\frac{\dd}{\dd y}\alpha \not\equiv 0$.

\begin{proof}

It is convenient to introduce an \emph{ad hoc} variable $s$ and to consider $x$ and $\alpha$ functions of 
$(F,s)$, where $G=s^2$ by definition. (This is possible because $\frac{\dd(F,G)}{\dd(x,\alpha)}(1,0)=32\neq 0$.)
Consider a differential form
$$\kappa=\frac{\dd x}{\dd s} dF+\left(\frac{\dd \alpha}{\dd s}-F\frac{\dd x}{\dd s}\right)\frac{ds}{s}.$$
It is analytic near the origin of the $(F,s)$ plane  (because $s^{-1}\frac{\dd}{\dd s}=2\frac{\dd}{\dd G}$). 
We have  the relation \eqref{e2a}  which in the new variables takes the form
$$F\frac{\dd x}{\dd F}+s\frac{\dd x}{\dd s}=\frac{\dd \alpha}{\dd F}.$$
It follows that
$$d\kappa=\frac{1}{s}\frac{\dd }{\dd s}\left(F\frac{\dd x}{\dd F}+s\frac{\dd x}{\dd s}-\frac{\dd \alpha}{\dd F}\right) ds\wedge dF=0.$$  
By the Poincare lemma, there is a unique analytic function $\Phi(F,s)$ such that $\Phi(0,0)=0$ and 
$\kappa=d\Phi$.

This  form is odd  with respect to the second variable,  in the sense that $\sigma^*\kappa=-\kappa$
where   $\sigma:(F,s)\mapsto (F,-s)$. If $\gamma$ is  a path connecting the origin $(0,0)$ to a given point $(F,s)$, then
$$\Phi(F,-s)=\int_{\sigma \gamma}\kappa=\int_{\gamma}\sigma^*\kappa=-\int_{\gamma}\kappa=-\Phi(F,s).$$
We have $\Phi^2(F,s)=\Phi^2(F,-s)$, hence  the  square  $\Phi^2$ is a well defined analytic function of $F$ and $G=s^2$.  Now we can change  the variables back and denote  $f(x,\alpha)=\Phi^2$. We have $f(1,0)=0$ (essentially, by  assumptions). 

Near the origin  $\Phi(F,s)=\left(\frac{1}{4}+O(F)\right)s+O(s^3)$,
hence 
$$\frac{\dd \Phi^2}{\dd F}(0,0)=0,\, \frac{\dd \Phi^2}{\dd G}(0,0)=\frac{1}{16},$$
and 
$$\frac{\dd}{\dd \alpha}f(1,0)=\left(\frac{3}{2}\frac{\dd}{\dd F}+8\frac{\dd}{\dd G}\right)\Phi^2(0,0)=\frac{1}{2}.$$
By the implicit function theorem, in a neighbourhood of the origin  there is a unique analytic  function of two  variables $\alpha(x,y)$ such that $\alpha(1,0)=0$ and
$$f(x, \alpha(x,y))=y^2.$$

From  \eqref{e2a} and the definition of  $\kappa$ we have 
$$sd\Phi=d\alpha-F dx.$$
Now, in the variables $(x,y)$ we  have $\Phi^2=y^2$, hence $d\Phi^2=dy^2$ and\footnote{As customary,  $dy^2$ actually  means $dy\otimes dy$, a tensor square.} 
$$(d\alpha-F(x,\alpha) dx)^2-G(x,\alpha)dy^2=0.$$
This  equality  implies \eqref{e2} (and is essentially  equivalent to  it).
\end{proof}

{\bf Remark 1}  In  Lemma 3, the condition $\alpha(1,0)=0$ is crucial. In this case  we cannot take the square root of the second equation, and  solving the system is more difficult then for $\alpha(1,0)>0$. Unfortunately, in the latter case the function $\alpha$ would have no extrema,  and the Euler flow $u$  in the following section could not be extended to the whole space.


{\bf Remark 2}  Note that  $\alpha(x,y)=\alpha(x,-y)$ pretty much by definition. An  interesting  consequence  is  $G(x,\alpha(x,0))=0$.  
(Which follows from  \eqref{e2} and $\tfrac{\dd }{\dd y}\alpha(x,0)=0$.) Because of this, the function  $\alpha(x,0)$ is actually another analytic  solution of \eqref{e2}.  


{\bf Remark 3}  In the given  proof,  the main technical difficulty is the absence of an  inverse  map  to $(F,s)\mapsto (x,\alpha)$. We circumvent this obstacle by  artificially constructing  a function $\Phi^2(F,s)$ with a well defined ``pullback''.     One of the referees  pointed out that there is a more straightforward  (although  not unrelated)  proof  using the Cartan - K\"{a}hler Theorem. We may consider a form $\omega=d\alpha- pdx-qdy$ on a manifold of dimension 3 defined by equations 
$$p=F(x,\alpha), q^2=G(x,\alpha)$$ 
in variables $x,y,\alpha, p, q$. This form satisfies the integrability condition $\omega\wedge d\omega=0$,   so we can use it to construct the  function
$\alpha(x,y)$.  (An important detail is that $\omega\slash q$ is analytic.)


{\bf Remark 4}   The method of the proof is  constructive and may be used to compute 
the Taylor series 
$$\alpha(x,y)=2(x-1)^2+2y^2+3(x-1)^3+3(x-1)y^2+O((|x-1|+|y|)^4).$$
The first two terms are important, so it may be appropriate to include a  direct computation of them.


$$$$
{\bf Lemma 4}
The function  $\alpha(x,y)$  has a strict  local minimum at  $(x,y)=(1,0)$.

\begin{proof}

The first  derivatives at this point are zero by \eqref{e2}. We have
$$\frac{\dd^2 }{\dd x^2}\alpha=\frac{\dd}{\dd x} F= \frac{\dd F}{\dd x}=4,\,\frac{\dd^2 }{\dd x \dd y}\alpha=
\frac{\dd}{\dd y} F=\left(\frac{\dd F}{\dd \alpha }\right)\left(\frac{\dd }{\dd y}\alpha\right)=0.$$
Finally, (at any  point) we have the equality
$$\left(\frac{\dd G}{\dd \alpha }\right)\left(\frac{\dd }{\dd y}\alpha\right)=\frac{\dd}{\dd y} G=2\left(\frac{\dd }{\dd y}\alpha\right)\left(\frac{\dd^2 }{\dd y^2}\alpha\right).$$
As the derivative $\frac{\dd }{\dd y}\alpha$ is not identically zero, it implies
$$\frac{\dd^2}{\dd y^2} \alpha=\frac{1}{2}\frac{\dd G}{\dd \alpha}.$$ 
At  the point under consideration we have then $\frac{\dd^2 }{\dd y^2}\alpha=4$.
The second differential $d^2\alpha$  is positively definite, so this point is a strict minimum.

\end{proof}

\section{The  flow}

We use the standard cylindrical coordinates. For a velocity field 
with axial symmetry the Euler equation $(u\cdot\nabla) u=-\nabla p$ takes   the form 

\begin{equation}
\left\{
\begin{array}{rl}
u_\rho\frac{\dd}{\dd\rho}u_\rho+u_z \frac{\dd}{\dd z}u_\rho-\frac{1}{\rho}u_\varphi^2=-\frac{\dd}{\dd\rho}p,\\
u_\rho\frac{\dd}{\dd\rho}u_\varphi+u_z \frac{\dd}{\dd z}u_\varphi+\frac{1}{\rho}u_\rho u_\varphi=0,\\
u_\rho\frac{\dd}{\dd\rho}u_z+u_z \frac{\dd}{\dd z}u_z=-\frac{\dd}{\dd z}p.
\end{array} \right. \label{e3}
\end{equation}

For  $R>0$, denote $a=\alpha\left(\frac{\rho}{R}, \frac{z}{R} \right)$. For the sake of convenience, we denote by $\CC$ the circle $\rho=R, z=0$ 
where $a=0$. Let 
\begin{equation}
p=\frac{aR^4}{4},\,b=\frac{R^3}{4}\sqrt{H(a)},\, u=\frac{1}{\rho}\left(\frac{\dd p}{\dd z}e_\rho-\frac{\dd p}{\dd \rho}e_z+be_\varphi\right).  \label{e4}
\end{equation}
(Note that $b$ is not smooth on $\CC$.) Obviously,  $\ddiv  u=0$ outside $\CC$.

The  fields $(u,p)$ given by  \eqref{e4}  satisfy \eqref{e3} in a neighbourhood of $\CC$  (but not on the curve itself).

\begin{proof}

The second equation of  \eqref{e3}  is obvious. The last one is equivalent to 
$$x\left(\frac{\dd }{\dd x}\alpha\right)\frac{\dd^2}{\dd x\dd y} \alpha+\left(\frac{\dd }{\dd y}\alpha\right)\left(\frac{\dd }{\dd x}\alpha+4x^3-x\frac{\dd^2 }{\dd x^2}\alpha\right)=0,$$
where $x=\frac{\rho}{R},\,y=\frac{z}{R}$ and $\alpha=\alpha(x,y)$. After multiplication  by  $\frac{\dd }{\dd y}\alpha$, using \eqref{e2} and
$$\frac{\dd G}{\dd x}+F\frac{\dd G}{\dd \alpha}=\frac{\dd }{\dd x}\left(\frac{\dd  }{\dd y}\alpha\right)^2=2\left(\frac{\dd }{\dd y}\alpha\right)\frac{\dd^2}{\dd x\dd y} \alpha$$
we have
$$\frac{1}{2}xF\left(\frac{\dd G}{\dd x}+F\frac{\dd G}{\dd \alpha}\right)+G\left(F+4x^3-x\left(\frac{\dd F}{\dd x}+F\frac{\dd F}{\dd \alpha}\right)\right)=0,$$
which follows from \eqref{e2a}, \eqref{e2b}.

Finally, the first equation is
$$x\left(\frac{\dd }{\dd x}\alpha\right)\frac{\dd^2}{\dd y^2} \alpha-x\left(\frac{\dd }{\dd y}\alpha\right)\frac{\dd^2}{\dd x\dd y} \alpha+\left(\frac{\dd }{\dd y}\alpha\right)^2-4x^3 \frac{\dd }{\dd x}\alpha+H(\alpha)=0,$$
or
$$\frac{1}{2}xF\frac{\dd G}{\dd \alpha}-\frac{1}{2}x\left(\frac{\dd G}{\dd x}+F\frac{\dd G}{\dd \alpha}\right)+G-4x^3 F+H(\alpha)=0,$$
which is again a consequence of \eqref{e2a}, \eqref{e2b}.
\end{proof}

As introduced, this  Euler flow is only defined in a vicinity of the circle $\CC$.  
However, this flow  satisfies an additional  condition $u\cdot\nabla p=0$  which is
very  useful for our purposes. Consider  another  field $\widetilde{u}=\omega(p)u$  where $\omega$ is a smooth function. Due to  the above condition 
we have 
$$\ddiv \widetilde{u}=\omega(p)\ddiv u+\omega'(p) (u\cdot\nabla p)=0$$
and 
$$(\widetilde{u}\cdot\nabla) \widetilde{u}=\omega^2(p)(u\cdot\nabla) u +\omega(p)\omega'(p)(u\cdot\nabla p)u=-\omega^2(p)\nabla p.$$
So, regardless of a choice of the  function $\omega$, the field  $\widetilde{u}$ is also an Euler flow, with the corresponding  pressure  determined by $d\widetilde{p}=\omega^2(p)\,dp$.

Due to Lemma 4, we may assume that $\omega=\omega(p)$  in a vicinity of the circle $\CC$ and $\omega=0$ outside  this domain.  
If  $\supp(\omega)\subset [\varepsilon, 2\varepsilon]$ (as a function of $p$)   with  $\varepsilon>0$  sufficiently small,  then we have $\widetilde{u}\in C^\infty(\Rs^3)$. The  new flow   is supported in a toroidal domain which can  be made arbitrarily close to the circle. This completes the proof of Theorem 1.

{\bf Remark 5} The poloidal stream function $\Psi=a$ is  a solution\footnote{This  is probably  what an expert would expect in this situation, but the author does not know an appropriate reference to make it a meaningful discussion.}
of the Grad-Shafranov equation in the following form  ($R=1$) 
$$(\dd_{\rho\rho}+\dd_{zz}-\frac{1}{\rho}\dd_\rho)\Psi=10\rho^2-\frac{1}{2}H'(\Psi).$$

\section{Generalized Beltrami  flows}

The condition $u\cdot\nabla p=0$  mentioned above  means that  the pressure   is 
constant along a streamline. This is very uncommon, and  the only other nontrivial example known to the author is 
a flow on the 3-sphere  constructed in \cite{KKP}. (The trivial examples are  vortices   (rotational flows) and their variations.)
By the Bernoulli theorem $|u|^2$ must also be a first integral; indeed, from  \eqref{e4}, \eqref{e2} and the definition of $G$  we have 
$$|u|^2=\frac{1}{\rho^2}\left[ \left(\frac{\dd p}{\dd z}\right)^2+\left(\frac{\dd p}{\dd \rho}\right)^2+b^2\right]=3p.$$
For the modified flow $\widetilde{u}$ the formula is different  but  $|\widetilde{u}|^2$ is still a function of the pressure $\widetilde{p}$.

Recall that for a  Beltrami  flow $u$  the Bernoulli function $B=p+\frac{1}{2}|u|^2$ is constant. 
The case when the Bernoulli function  depends  on the pressure  may be considered a  generalization, and
constructed flows  belongs to this category  ($B=\frac{5}{2}p$ for the original flow). As was pointed out by Arnold \cite{A}[II.1.B], for a non-constant $B$ 
both the streamlines  and the vortex  lines lie on the surfaces $B=\text{const}$; in our situation these are the same 
as $p=\text{const}$. It  makes a difference because  in this case the  flow  sheets   become 
independent   in a  sense,  so the flow  may be ``modulated''  (a trick we used in the previous section).

One of the referees pointed out to  the author  that  a generalized Beltrami  flow (with an extremum of pressure at some point)  has a peculiar restriction on the behaviour  of the pressure.  Let  $(u,p)$ be  such a flow, and  assume that $|u|^2=3p$ as before (we can do  this  without loss of generality).  By the same recipe as above  we may then construct another flow $(\widetilde{u},\widetilde{p})$.   
If it has compact support then (e.g. \cite{CC})
$$\int_{\Rs^3} (|\widetilde{u}|^2+3\widetilde{p})\,d x=0.$$
To make sense of this it is convenient to introduce a function  $V(c)=Vol(\{x\in \Rs^3: p( x)\le c\})$. The equality then  becomes
$$\int_{\Rs^3} (p\omega^2(p)+\widetilde{p})\,d  x=\int \omega^2(p) (p\,dV(p) -V(p)\,dp)=0.$$
It must be frue for any  function $\omega$ which means  $V(p)=p\cdot\text{const}$.

\end{document}